\documentclass[12pt]{article}

\usepackage{latexsym}
\usepackage{amsfonts}
\usepackage{amssymb}
\usepackage{amsmath}
\usepackage{mathrsfs}
\usepackage{hyperref}

\newcommand{\be}{\begin{equation}}
\newcommand{\ee}{\end{equation}}
\newcommand{\bea}{\begin{eqnarray}}
\newcommand{\eea}{\end{eqnarray}}
\newcommand{\bean}{\begin{eqnarray*}}
\newcommand{\eean}{\end{eqnarray*}}
\newcommand{\brray}{\begin{array}}
\newcommand{\erray}{\end{array}}

\newtheorem{dfn}{Definition}[section]
\newtheorem{thm}[dfn]{Theorem}
\newtheorem{lmma}[dfn]{Lemma}
\newtheorem{ppsn}[dfn]{Proposition}
\newtheorem{crlre}[dfn]{Corollary}
\newtheorem{xmpl}[dfn]{Example}
\newtheorem{rmrk}[dfn]{Remark}

\newcommand{\bdfn}{\begin{dfn}\rm}
\newcommand{\bthm}{\begin{thm}}
\newcommand{\blmma}{\begin{lmma}}
\newcommand{\bppsn}{\begin{ppsn}}
\newcommand{\bcrlre}{\begin{crlre}}
\newcommand{\bxmpl}{\begin{xmpl}}
\newcommand{\brmrk}{\begin{rmrk}\rm}

\newcommand{\edfn}{\end{dfn}}
\newcommand{\ethm}{\end{thm}}
\newcommand{\elmma}{\end{lmma}}
\newcommand{\eppsn}{\end{ppsn}}
\newcommand{\ecrlre}{\end{crlre}}
\newcommand{\exmpl}{\end{xmpl}}
\newcommand{\ermrk}{\end{rmrk}}

\author{S. Sundar}
\title{On the Wiener-Hopf compactification of a symmetric cone}

\begin{document}
\maketitle 
\begin{abstract}

Let $V$ be a finite dimensional real Euclidean Jordan algebra with the identity element $1$. Let $Q$ be the closed convex cone of squares.  We show that the Wiener-Hopf compactification of $Q$ is the interval $\{x \in V: -1 \leq x \leq 1\}$. As a consequence, we deduce that the $K$-groups of the Wiener-Hopf $C^{*}$-algebra associated to $Q$ are trivial.
\end{abstract}
\noindent {\bf AMS Classification No. :} {Primary 46L80; Secondary 17CXX.}  
\section{Introduction}
Let $P$ be a closed convex  cone of $\mathbb{R}^{n}$. Assume that $P$ is spanning i.e. $P-P=\mathbb{R}^{n}$. For $f \in C_{c}(\mathbb{R}^{n})$, let $W(f)$ be the operator on $L^{2}(P)$ defined by
\[
W(f)\xi(s)=\int_{t \in P}f(s-t)\xi(t) dt
\]
for $\xi \in L^{2}(P)$. The $C^{*}$-algebra generated by $\{W(f):f \in C_{c}(\mathbb{R}^{n})\}$ is called the Wiener-Hopf $C^{*}$-algebra associated to the cone $P$ and is denoted $\mathcal{W}(P)$. The study of $\mathcal{W}(P)$ from the groupoid perspective 
was initiated by Muhly and Renault in \cite{Renault_Muhly} and further developed by Nica in \cite{Nica_WienerHopf}.  

Let us recall the groupoid model for $\mathcal{W}(P)$ obtained in \cite{Renault_Muhly}. Let $Y$ be the unit ball of $
L^{\infty}(\mathbb{R}^{n})$ endowed with the weak $*$-topology. The group $\mathbb{R}^{n}$ acts on $Y$ by translations. Let $X$ be the closure of $\{1_{a-P}: a \in P\}$ in $Y$.  Here as usual, for  $A\subset \mathbb{R}^{n}$, $1_{A}$ denotes the characteristic function 
of $A$. Note that $X$ is compact and is left invariant by the additive semigroup $P$.  We call $X$ as the Wiener-Hopf compactification of $P$.  Let $\mathcal{G}$ be the reduction of the transformation groupoid $Y \rtimes \mathbb{R}^{n}$ onto $X$. 
Then it is proved in \cite{Renault_Muhly} that $\mathcal{W}(P)$ is isomorphic to the reduced $C^{*}$-algebra of the groupoid $\mathcal{G}$. For a self-contained treatment of these results for a general Ore semigroup, we refer the reader to \cite{Jean_Sundar}.
For a generalization in the context of ordered homogeneous spaces, see \cite{Hilgert_Neeb}.

Two main classes of cones studied in \cite{Renault_Muhly} are polyhedral cones and self-dual homogeneous cones, also called as symmetric cones. In particular, their Wiener-Hopf compactifications are described in detail.
The groupoid model for $\mathcal{W}(P)$, in the case of polyhedral cones, is exploited
to obtain index theorems in \cite{AJ07}, \cite{AJ08}. Moreover in \cite{All11}, it is proved that the $K$-groups of $\mathcal{W}(P)$ vanish when $P$ is polyhedral. Here we prove that for a self-dual homogeneous cone the $K$-groups of the 
Wiener-Hopf $C^{*}$-algebra vanish. This requires us to take another look at the Wiener-Hopf compactification of a symmetric cone. 

Let $V$ be a finite dimensional real Euclidean Jordan algebra with identity $1$. Denote the closed  convex cone of squares by $Q$.  Let us recall the description of the Wiener-Hopf compactification of Q given in \cite{Renault_Muhly}. For an idempotent $e$ in $V$, let $e^{\perp}=1-e$. 
Let $P:V \to End(V)$ be the quadratic map which determines the structure of the Jordan algebra.  Let $Y:=\{(e,x) \in V \times Q: e^{2}=e, e^{\perp}x=x\}$. For $(e,x) \in Y$, let \[A_{(e,x)}=\{v \in V: x+P(e^{\perp})v \geq 0\}.\] 
Then it is proved in \cite{Renault_Muhly} that the map $Y \ni (e,x) \to -A_{(e,x)} \in X$ is a bijection. Also note that the map $Y \ni (e,x) \to e+(x-e^{\perp})(x+e^{\perp})^{-1} \in V$, (where the expression $(x+e^{\perp})^{-1}$ stands for the inverse of $x+e^{\perp}$ in the Jordan subalgebra $P(e^{\perp})V$), is injective with the interval $[-1,1]:=\{u \in V: -1 \leq u \leq 1\}$ as the range. It is then a natural question to ask whether the topology on $[-1,1]$, transported from the weak $*$-topology on $X$, is the subspace topology on $[-1,1]$.  The main objective of this short paper is to answer this in the affirmative. 

We should remark that  this could  probably be proved directly following the ideas in \cite{AJ07} (Theorem 9).  However, we follow a different path by appealing to the axiomatic characterisation of the Wiener-Hopf compactification of an Ore semigroup (Prop 5.1, \cite{Jean_Sundar}). We believe that this axiomatic approach is elegant and has potential applications to semigroups other than convex cones. We illustrate this by proving that the Wiener-Hopf algebra associated to the continuous "ax+b"-semigroup has vanishing $K$-theory.

We should also mention here that for special Jordan algebras, i.e. Jordan algebras that admit an embedding into the Jordan algebra of self-adjoint matrices, the Wiener-Hopf compactification of the positive elements is identified, via the Cayley transform, with a subset of the unitary group in \cite{Sundar_Toeplitz}. The same could be done here by passing to the complexification. Since we are only interested in compactifying $Q$, we avoid passing to the complexification of $V$ and instead work with the transform $Q \ni x \to \frac{x-1}{x+1} \in V$.

\section{Preliminaries}

Let $P$ be a topological semigroup and $X$ be a topological space. By an action of $P$ on $X$, we mean a continuous map $X \times P \to X$, the image of $(x,a)$ is denoted $xa$, such that $x(ab)=(xa)b$ for $x \in X$ and $a,b \in P$. If $P$ has an identity element $e$ , we demand that $xe=x$ for $x \in X$. We say that the action is injective if for each $a \in P$, the map $X \ni x \to xa \in X$ is injective.  

We reserve the letter $G$ to denote a locally compact, Hausdorff, second countable topological group. Let $P\subset G$ be a closed subsemigroup. We say that 
\begin{enumerate}
\item[(1)] $P$ is solid if the interior of $P$, denoted $Int(P)$, is dense in $P$, and
\item[(2)] $P$ is right Ore if $PP^{-1}=G$.
\end{enumerate}  
We start with an elementary lemma.

\begin{lmma}
\label{extension}
Let $P \subset G$ be a closed, solid subsemigroup of $G$ containing the identity element $e$ of $G$. 
Let $X$ be a compact metric space and let $\tau:X \times Int(P) \to X$ be an injective action. Then there exists a unique injective action $\widetilde{\tau}:X \times P \to X$ such that $\widetilde{\tau}(x,a)=\tau(x,a)$ for $x \in X$ and $a \in Int(P)$.  
\end{lmma}
\textit{Proof.} As usual, we write $\tau(x,a)$ as $xa$ for $x \in X$ and $a \in Int(P)$. Let $x \in X$ and $a \in P$ be given. Let $(a_{n})$ be a sequence in $Int(P)$ such that $a_{n}\to a$. We claim that $(xa_{n})$ converges. Since $X$ is a compact metric space, it is enough to show that every convergent subsequence of $(xa_{n})$ converges to the same limit. Thus suppose $(xa_{n_k})$ and $(xa_{m_k})$ be subsequences of $(xa_{n})$ converging to $y$ and $z$ respectively. Let $s \in Int(P)$ be given. Note that $P Int(P) \subset Int(P)$. Since $a_{n}s \to as$ and the action of $Int(P)$ on $X$ is continuous, it follows that $ys=x(as)=zs$. Since the map $X \ni x \to xs \in X$ is injective, it follows that $y=z$. Hence the sequence $(xa_{n})$ converges. 

Now let $(b_{n})$ and $(c_{n})$ be sequences in $Int(P)$ such that $b_{n} \to a$ and $c_{n} \to a$. Suppose that $(xb_{n}) \to y$ and $(xc_{n}) \to z$.  Let $s \in Int(P)$ be given.  Since $b_{n}s \to as$ and $c_{n}s \to as$, it follows, by the continuity of the action of $Int(P)$, that $ys=x(as)=zs$. Hence $y=z$. 

Define for $x \in X$ and $a \in P$, $\widetilde{\tau}(x,a)= \lim_{n \to \infty}\tau(x,a_{n})$ where $(a_{n})$ is any sequence in $Int(P)$ converging to $a$. We have shown that $\widetilde{\tau}$ is well-defined. We leave it to the reader to verify that $\widetilde{\tau}$ is an injective action of $P$ on $X$ and it extends $\tau$.  Uniqueness of the action follows from the fact that $Int(P)$ is dense in $P$. \hfill $\Box$. 

Let $P \subset G$ be a closed subsemigroup containing the identity element. Assume that $P$ is solid and right Ore. We recall the Wiener-Hopf compactification associated to the pair $(P,G)$.
 Let $X$ be a compact Hausdorff space on which $P$ acts injectively. For $x \in X$, let \[A_{x}:=\{g \in G: \textrm{there exists $a,b \in P$ and $y \in X$ such that $g=ab^{-1}$ and $xa=yb$}\}.\] Note that since $PInt(P) \subset Int(P)$, for $x \in X$, $A_{x}$ can alternatively be described as
\[
A_{x}:=\{g \in G: \textrm{there exists $a,b \in Int(P)$, $y \in X$ such that $g=ab^{-1}$ and $xa=yb$}\}.\]

\begin{dfn}
\label{Wiener-Hopf compactification}
Let $X$ be a compact Hausdorff space on which $P$ acts injectively. The space $X$ is called the 
Wiener-Hopf compactification of $P$ if the following conditions are satisfied.
\begin{enumerate}
\item[(C1)] For $a \in P$, $\{x(ba): b \in Int(P), x \in X\}$ is open in X. 
\item[(C2)] There exists $x_{0} \in X$ such that $A_{x_0}=P$ and $\{x_{0}a: a \in P\}$ is dense in $X$.
\item[(C3)] For $x,y \in X$, if $A_{x}=A_{y}$ then $x=y$.
\end{enumerate}
\end{dfn}
\begin{rmrk}
An explicit model for the Wiener-Hopf compactification of the pair $(P,G)$ is constructed in \cite{Renault_Muhly}. The closure of $\{1_{P^{-1}a}: a \in P\}$ in $L^{\infty}(G))$ is a model for the Wiener-Hopf compactification. Here $L^{\infty}(G)$ is given the weak$^{*}$-topology and the action of $P$ is by right translation. 
For a self-contained proof of  the existence  and the uniqueness of the Wiener-Hopf compactification, we refer the reader to  the article \cite{Jean_Sundar}, in particular to Prop.5.1, \cite{Jean_Sundar}. 
The term "Order compactification" is used instead of the Wiener-Hopf compactification in \cite{Jean_Sundar}.
\end{rmrk}

\begin{rmrk}
Prop 5.1 of \cite{Jean_Sundar} requires that the map $X \times Int(P) \ni (x,a) \to xa \in X$ is open. However it is proved in Theorem 4.3 of \cite{Jean_Sundar} that the openness of the map $X \times Int(P) \ni (x,a) \to xa \in X$ is equivalent to (C1).
\end{rmrk}

We recall the essentials of Jordan algebras and fix notations that is needed to read this paper. We refer the reader to the monograph \cite{Faraut} for proofs. 

Let $V$ be a finite dimensional real Euclidean Jordan algebra with the identity element $1$. For $x \in V$, let $L(x):V \to V$ be defined by $L(x)y=xy$. Let $P:V \times V \to End(V)$ be the  bilinear map defined by the equation $P(x,y):=L(x)L(y)+L(y)L(x)-L(xy)$. For $x \in V$, we simply denote $P(x,x)$ by $P(x)$.  

 Let $u \in V$ be invertible. For $a,b \in V$, let $a\star_{u}b=P(a,b)u$. Then $\star_{u}$ is a Jordan product on $V$. We denote the Jordan algebra $(V,\displaystyle\star_{u})$ by $V_{u}$. The Jordan algebra $V_{u}$ is called the mutation at $u$. The Jordan algebra $V_{u}$ is unital with identity $u^{-1}$. Moreover, for $x \in V$, $x$ is invertible in $V$ if and only if $x$ is invertible in $V_{u}$. For an invertible element $x \in V$, we write the inverse of $x$ in $V_{u}$ by $x_{u}^{-1}$. Then $x_{u}^{-1}=P(u^{-1})x^{-1}$.   
 We also need Hua's identity. Suppose $a,b,a+b \in V$ are invertible, then $a+P(a)b^{-1}$ is invertible and 
 \[
 (a+b)^{-1} + (a+P(a)b^{-1})^{-1}=a^{-1}.
 \]

Assume that $V$ is Euclidean i.e. there exists an inner product $\langle,\rangle$ with respect to which $L(x)$ is symmetric for every $x \in V$.  Let $Q:=\{x^{2}:x \in V\}$. It is well known that $Q=\{x \in V: L(x) \textrm{~is positive}\}$ and that $Q$ is a closed convex cone. Denote the interior of $Q$ by $\Omega$. Then $\Omega=\{x \in Q: x \textrm{~is invertible}\}$.  Let $x,y \in V$ be given. We write $x \leq y$  if $y-x \in Q$ and write $x<y$ if $y-x \in \Omega$. 

For $x \in V$, let \[\sigma(x):=\{\lambda \in \mathbb{R}: x-\lambda.1 ~\textrm{is not invertible}\}.\] 
The set $\sigma(x)$ is called the spectrum of $x$.  An element $e \in V$ is called an idempotent if $e^{2}=e$. Let $e,f \in V$ be idempotents. We say that $e$ and $f$ are orthogonal if $ef=0$.   The spectral theorem holds in finite dimensional Euclidean Jordan algebras i.e. given $x \in V$, there exist orthogonal idempotents, necessarily unique, $\{c_{\lambda}: \lambda \in \sigma(x)\}$ such that $x=\displaystyle \sum_{\lambda \in \sigma(x)}\lambda c_{\lambda}$ and $\displaystyle \sum_{\lambda \in \sigma(x)}c_{\lambda}=1$.

We recall now the Peirce decomposition. Let $e \in V$ be an idempotent.  For $\lambda \in \mathbb{R}$, let $V_{\lambda}(e):=\{x \in V: ex=\lambda x\}$. Then $V$ decomposes orthogonally as \[V=V_{0}(e)\oplus V_{1}(e) \oplus V_{\frac{1}{2}}(e)\] and the decomposition is called the Peirce decomposition of $V$ w.r.t. the idempotent $e$.  The orthogonal projection onto $V_{1}(e)$ is given by $P(e)$.  Also $V_{1}(e)$ is a Jordan subalgebra of $V$ and $e$ is the identity of $V_{1}(e)$. Let $x \in V_{1}(e)$ be given.  If $x$ is invertible in $V_{1}(e)$, we denote the inverse of $x$ in $V_{1}(e)$ as $_{e}x^{-1}$. For an idempotent $e$ in $V$, set $e^{\perp}:=1-e$. Note that $V_{0}(e^{\perp})=V_{1}(e)$.

\section{The Wiener-Hopf Compactification of a symmetric cone}
For the rest of this paper, let $V$ be a finite dimensional real Euclidean Jordan algebra with identity $1$ and denote the closed convex cone of squares by $Q$. Denote the interior of $Q$ by $\Omega$. 
Let \[X:=\{u \in V: -1 \leq u \leq 1\}\]and endow $X$ with the subspace topology inherited from the topology on $V$. Then $X$ is compact. We also write $[-1,1]$ to denote the set $X$.   
We show here that $X$ is the Wiener-Hopf compactification of the cone $Q$.
First we describe the action of $Q$ on $X$. In view of Lemma \ref{extension}, it is enough to describe the action of $\Omega$ on $X$.  For $u \in X$ and $a \in \Omega$, let 
\[
u \boxplus a:= 1-2\widetilde{u}+2P(\widetilde{u})(\widetilde{u}+a^{-1})^{-1}
\]
where $\widetilde{u}:=\frac{1-u}{2}$.  

\begin{rmrk}
Note that for $u \in X$, $\widetilde{u} \geq 0$. Also observe that if $a \in Q$ and $b \in \Omega$ then $a+b \in \Omega$. Hence the map $\boxplus:X \times \Omega \to V$ is 
well-defined. Since the inversion is continuous, it follows that the map $\boxplus:X \times \Omega \to V$ is continuous.
\end{rmrk}

Let $i:Q \to X$ be defined by $i(x)=\frac{x-1}{x+1}$. 
Note that $i$ is an embedding, $i(Q)$ is dense in $X$ and $i(Q):=\{u \in V: -1 \leq u <1\}$.  These are straightforward 
consequences of the spectral theorem and thus we leave the verifications of these statements to the reader. 

\begin{ppsn}
\label{action}
For $u \in X$ and $a \in \Omega$, $u \boxplus a \in X$. Also the map $\boxplus:X \times \Omega \to X$ is an  action of $\Omega$ on $X$.
\end{ppsn}
\textit{Proof.} Let $x \in Q$ and $a \in \Omega$ be given. We claim that $i(x)\boxplus a=i(x+a)$.  Let $u:=i(x)=1-\frac{2}{x+1}$. Then $\widetilde{u}=(x+1)^{-1}$ i.e. $x+1=\widetilde{u}^{-1}$. 
Now note that 
\begin{align*}
 i(x+a)&=1-\frac{2}{x+a+1} \\
 &=1-2(\widetilde{u}^{-1}+a)^{-1} \\ 
 &=1-2\Big(\widetilde{u}-\big(\widetilde{u}^{-1}+P(\widetilde{u}^{-1})a^{-1}\big)^{-1}\Big) ~~\Big(~\textrm{by Hua's identity}\Big) \\
 &=1-2\widetilde{u}+2\big(P(\widetilde{u})^{-1}\widetilde{u}+P(\widetilde{u})^{-1}a^{-1}\big)^{-1} \\
 &=1-2\widetilde{u}+2P(\widetilde{u})(\widetilde{u}+a^{-1})^{-1}.
 \end{align*}
Hence $i(x)\boxplus a=i(x+a)$.  This proves the claim. 

Since $i(Q)$ is dense in $X$, $X$ is closed in $V$ and the map $\boxplus:X \times \Omega \to V$ is continuous, it follows that 
$u\boxplus a \in X$ for $u \in X$ and $a \in \Omega$.  The equality $i(x)\boxplus a=i(x+a)$  implies that 
$(u \boxplus a)\boxplus b=u \boxplus(a+b)$ for $u \in i(Q)$ and $a,b \in \Omega$.   Now the equality $(u \boxplus a)\boxplus b =u \boxplus (a+b)$ for $u \in X$ and $a,b \in \Omega$ 
follows from the continuity of the map $\boxplus$ and the density of $i(Q)$ in $X$.  This completes the proof. \hfill $\Box$

The fact that $\boxplus$ is an injective action requires a bit more work. We now proceed towards proving it. We need the well known fact that there exists an open set $U \subset V$ such that $0 \in U$ and
if $x \in U$, $1+x$ is invertible and 
\[
(1+x)^{-1}=1-x+x^{2}-x^{3}+\cdots
\]
The following lemma is crucial in what follows.
\begin{lmma}
\label{crucial}
Let $e \in V$ be an idempotent and let $V=V_{1}(e)\oplus V_{0}(e) \oplus V_{\frac{1}{2}}(e)$ be the Peirce decomposition w.r.t. the idempotent $e$. Let $a \in \Omega$ be given.  Set $a_{0}:=P(e)a$. 
Then 
\begin{enumerate}
\item[(1)] $a_{0}$ is invertible in the Jordan algebra $V_{1}(e)$, and
\item[(2)] $P(e)(u+a^{-1})^{-1}=~ _e(u+~_{e}a_{0}^{-1})^{-1}$ for $u \in V_{1}(e) \cap Q$.
\end{enumerate}

\end{lmma}
\textit{Proof.} Recall that if $x \in V_{1}(e)$ is invertible in $V_{1}(e)$, we denote the inverse of $x$ in $V_{1}(e)$ by $_{e}x^{-1}$.  Denote the quadratic representation of the Jordan algebra $V_{1}(e)$ on $V_{1}(e)$ by
$P_{e}$. Then $P_{e}(x)=P(e)P(x)P(e)$ for $x \in V_{1}(e)$. Since $P(e)$ is a projection, it follows that \[
P_{e}(a_{0})=P(e)P(P(e)a)P(e)=P(e)P(e)P(a)P(e)P(e)=P(e)P(a)P(e).\]
Since $a \in \Omega$, $P(a)$ is invertible and $P(a)$ is positive. Hence the cutdown $P(e)P(a)P(e)$ is invertible on $V_{1}(e)$. Thus $P_{e}(a_{0})$ is invertible and consequently $a_{0}$ is invertible in $V_{1}(e)$. 
This proves $(1)$.

Let \[A:=\{u \in V_{1}(e) \cap Q: P(e)(u+a^{-1})^{-1}=~_{e}(u+~_{e}a_{0}^{-1})^{-1} \}.\] Note that $0 \in A$. We claim that $A$ is a clopen subset of $V_{1}(e) \cap Q$.  Since the inversion is continuous, it follows that $A$ is 
closed.   

Let $u_{0} \in A$ be given. Set $w_{0}:=(u_{0}+a^{-1})^{-1}$ and $w_{1}:=~_{e}(u_{0}+~_{e}a_{0}^{-1})^{-1}$.  Since $u_{0} \in A$, $P(e)w_{0}=w_{1}$.
Let $V_{w_0}$ be the mutation of $V$ at $w_{0}$.  Let $U_{1} \subset V$ be an open subset containing $0$ such that if $u \in U_{1}$ then $u+w_{0}^{-1}$ is invertible in $V_{w_0}$ and 
\[
(u+w_{0}^{-1})_{w_{0}}^{-1}=w_{0}^{-1}-u + u \star_{w_{0}}u - u\star_{w_{0}}u \star_{w_0} u + \cdots.
\]
Thus for $u \in U_{1}$, $u+w_{0}^{-1}$ is invertible in $V$ and
\begin{align*}
(u+w_{0}^{-1})^{-1}&=P(w_{0})(u+w_{0}^{-1})_{w_0}^{-1}\\
                              &=w_{0}-P(w_{0})u+P(w_{0})(u \star_{w_{0}}u)-P(w_{0})(u \star_{w_0} u \star_{w_{0}} u)+\cdots.
\end{align*}

Working in the mutation  of $V_{1}(e)$ at $w_{1}$, we see that there exists an open subset  $U_{2}$ of $V$ containing $0$ such that if $u \in U_{2} \cap V_{1}(e)$ then $u+ _{e}w_{1}^{-1}$ is invertible in $V_{1}(e)$ and
\[
_{e}(u+~_{e}w_{1}^{-1})^{-1}=w_{1}-P(w_{1})u+P(w_{1})(u \star_{w_1}u)-P(w_{1})(u \star_{w_{1}}u \star_{w_1} u)+ \cdots.
\]
Observe that for $u,v \in V_{1}(e)$,  
\begin{align*}
u\star_{w_0}v&=P(u,v)w_{0} \\
                      &=P(P(e)u,P(e)v)w_{0} \\
                      &=P(e)P(u,v)P(e)w_{0} \\
                      &=P(e)P(u,v)w_{1} \\
                      &=P(e)(u\star_{w_{1}}v).
\end{align*}
Let $U:=U_{1} \cap U_{2} \cap V_{1}(e) \cap Q$.  Then $U$ is an open subset of $V_{1}(e) \cap Q$ containing $0$. Consider an element $u \in U$. Calculate as follows to observe that
\begin{align*}
&P(e)(u+u_{0}+a^{-1})^{-1} \\
&= P(e)(u+w_{0}^{-1})^{-1} \\
                                         &=P(e)(w_{0}-P(w_{0})u+P(w_{0})(u \star_{w_{0}}u)-P(w_{0})(u \star_{w_{0}} u \star_{w_{0}}u)+ \cdots)\\
                                         &=P(e)w_{0}-P(e)P(w_{0})P(e)u+P(e)P(w_{0})P(e)(u\star_{w_{1}}u)-P(e)P(w_{0})P(e)(u \star_{w_{1}}u \star_{w_1} u)+ \cdots \\
                                         &=P(e)w_{0}-P(P(e)w_{0})u+P(P(e)w_{0})(u\star_{w_1}u)-P(P(e)w_{0})(u\star_{w_1}u\star_{w_1}u)+ \cdots \\
                                         &=w_{1}-P(w_{1})u+P(w_{1})(u \star_{w_1}u)-P(w_{1})(u \star_{w_{1}}u \star_{w_1} u)+ \cdots \\
                                         &=~_{e}(u+~_{e}w_{1}^{-1})^{-1} \\
                                         &=~_{e}(u+u_{0}+~_{e}a_{0}^{-1})^{-1} 
       \end{align*}
 This implies that $u_{0}+U \subset A$. As as result, it follows that $A$ is open. Hence $A$ is a clopen subset of $V_{1}(e) \cap Q$. 
 But note that $V_{1}(e) \cap Q$ is convex and hence connected.  As a consequence, we 
 deduce that $A=V_{1}(e) \cap Q$.  This completes the proof. \hfill $\Box$

Let $Y:=\{(e,x) \in V \times V: e^{2}=e, x \in V_{0}(e) \cap Q\}$. The compact set $X$ can be put in bijective correspondence with $Y$ as follows. 
For an idempotent $e \in V$, let $e^{\perp}=1-e$.  Let $i:Y \to X$ be defined as
\[
i(e,x):=e+~_{e^{\perp}}(x+e^{\perp})^{-1}(x-e^{\perp}).
\]
For $x \in Q$, $i(x)$ stands for the case when $e=0$, i.e. $i(x)=i(0,x)=\frac{x-1}{x+1}$.  

Let $(e,x) \in Y$ be given. Let $x:=\sum_{i=1}^{n}\lambda_{i}f_{i}$ be the spectral decomposition of $x$ in $V_{1}(e)$ i.e. $\lambda_{i}$'s  are distinct, $f_{i}$'s are non-zero idempotents which are mutually orthogonal and $\sum_{i=1}^{n}f_{i}=e^{\perp}$.
Then $i(e,x)=e+\sum_{i=1}^{n}\frac{\lambda_{i}-1}{\lambda_{i}+1}f_{i}$.  Now it is clear that $i(e,x) \in X$.  Also $e$ is the spectral projection of $i(e,x)$ corresponding to the "eigenvalue" $1$. 

Let $(e,x), (f,y) \in Y$ be given.  Suppose $i(e,x)=i(f,y)$. Since the spectral projections of $i(e,x)$ and $i(f,y)$   corresponding to the eigenvalue $1$ are $e$ and $f$ respectively, it follows that $e=f$. Now the equation $i(e,x)=i(e,y)$ implies that $_{e^{\perp}}(x+e^{\perp})^{-1}(x-e^{\perp})=~_{e^{\perp}}(y+e^{\perp})^{-1}(y-e^{\perp})$ which in turn implies that $x=y$. This proves that $i$ is injective. 

Let $u \in X$ be given. Since $-1\leq u \leq 1$, the spectrum of $u$ is contained in $[-1,1]$. Let $e$ be the spectral projection of $u$ corresponding to the eigenvalue $1$.  Then there exits non-zero orthogonal idempotents $f_{1},f_{2},\cdots f_{n}$ and distinct scalars $\mu_{1},\mu_{2},\cdots ,\mu_{n} \in [-1,1)$ such that $u=e+\sum_{i=1}^{n}\mu_{i}f_{i}$ and $\sum_{i=1}^{n}f_{i}=e^{\perp}$. Let $\lambda_{i} \geq 0$ be such that $\mu_{i}=\frac{\lambda_{i}-1}{\lambda_{i}+1}$. Set $x:=\sum_{i=1}^{n}\lambda_{i}f_{i}$. Then $x \in V_{1}(e^{\perp})=V_{0}(e)$.  Clearly $i(e,x)=u$. This proves that $i$ is surjective.

\begin{lmma}
\label{explicit}
Let $(e,x) \in Y$ and $a \in \Omega$ be given. Then \[i(e,x)\boxplus a=i(e,x+P(e^{\perp})a).\] 

\end{lmma}
\textit{Proof.} We can assume that $e \neq 0$. For the case $e=0$ is already proved in Proposition \ref{action}. Let $u:=i(e,x)$ and set $v:=~_{e^{\perp}}(x+e^{\perp})^{-1}(x-e^{\perp})$. Then $u=e+v$ and $v \in V_{1}(e^{\perp})$. 
Let $\widetilde{u}:=\frac{1-u}{2}$ and $\widetilde{v}:=\frac{e^{\perp}-v}{2}$. Then $\widetilde{u}=\widetilde{v}$ and $\widetilde{v}=~_{e^{\perp}}(x+e^{\perp})^{-1}$.  Set $a_{0}=P(e^{\perp})a$. 

Recall that the projection onto $V_{1}(e^{\perp})$ is $P(e^{\perp})$ and
the one  onto $V_{1}(e)$ is $P(e)$. Also the spaces $V_{1}(e^{\perp})$ and $V_{1}(e)$ are orthogonal i.e. $P(e)P(e^{\perp})=0$. Hence $P(e)\widetilde{u}=0$ and $P(e^{\perp})\widetilde{u}=\widetilde{u}$.
Observe that 
\begin{align*}
P(e)(u\boxplus a)&=P(e)(1-2\widetilde{u}+2P(\widetilde{u})(\widetilde{u}+a^{-1})^{-1}) \\
                            &=e-2P(e)\widetilde{u}+2P(e)P(P(e^{\perp})\widetilde{u})(\widetilde{u}+a^{-1})^{-1} \\
                            &=e-0+2P(e)P(e^{\perp})P(\widetilde{u})P(e^{\perp})(\widetilde{u}+a^{-1})^{-1} \\
                            &=e.
\end{align*}
Also 
\begin{align*}
P(e^{\perp})(u \boxplus a)&=e^{\perp}-2P(e^{\perp})\widetilde{u}+2P(e^{\perp})P(P(e^{\perp})\widetilde{u})(\widetilde{u}+a^{-1})^{-1} \\
                                         &=e^{\perp}-2\widetilde{u}+2P(e^{\perp})P(e^{\perp})P(\widetilde{u})P(e^{\perp})(\widetilde{u}+a^{-1})^{-1} \\
                                          &=e^{\perp}-2\widetilde{u}+2P(e^{\perp})P(\widetilde{u})P(e^{\perp})(\widetilde{u}+a^{-1})^{-1} \\
                                          &=e^{\perp}-2\widetilde{u}+2P(P(e^{\perp})\widetilde{u})(\widetilde{u}+a^{-1})^{-1} \\
                                          &=e^{\perp}-2\widetilde{u}+2P(\widetilde{u})(\widetilde{u}+a^{-1})^{-1 }.
                                          \end{align*}
   Thus $u\boxplus a=P(e)(u\boxplus a)+P(e^{\perp})(u\boxplus a)$.  In other words, $u \boxplus a \in V_{1}(e)\oplus V_{0}(e)$. Now calculate by starting from the third equality of the previous calculation to observe that
   \begin{align*}
   P(e^{\perp})(u \boxplus a)&=e^{\perp}-2\widetilde{u}+2P(e^{\perp})P(\widetilde{u})P(e^{\perp})(\widetilde{u}+a^{-1})^{-1} \\
                                             &= e^{\perp}-2\widetilde{v}+2P(P(e^{\perp})\widetilde{u})P(e^{\perp})(\widetilde{v}+a^{-1})^{-1} \\
                                             &=e^{\perp}-2\widetilde{v}+2P(\widetilde{v})~_{e^{\perp}}(\widetilde{v}+~_{e^{\perp}}a_{0}^{-1})^{-1} ~~(\textrm{ by Lemma \ref{crucial}}).   
   \end{align*}                                                                               
Note that $\widetilde{v}$ is invertible in $V_{1}(e^{\perp})$.  (To avoid obscure notations, we will simply write $x^{-1}$ in place of $_{e}x^{-1}$ in the following calculation if $x \in V_{1}(e^{\perp})$ is invertible. )
Now calculate in the Jordan algebra $V_{1}(e^{\perp})$ to observe that
\begin{align*}
P(e^{\perp})(u \boxplus a)&=e^{\perp}-2\widetilde{v}+2P(\widetilde{v})(\widetilde{v}+a_{0}^{-1})^{-1} \\
                                          &=e^{\perp}-2\widetilde{v}+2(\widetilde{v}^{-1}+P(\widetilde{v}^{-1})a_{0}^{-1})^{-1} \\
                                          &=e^{\perp}-2\widetilde{v}+2(\widetilde{v}-(\widetilde{v}^{-1}+a_{0})^{-1})~~(\textrm{by Hua's identity in $V_{1}(e^{\perp})$}) \\
                                          &=e^{\perp}-2(x+e^{\perp}+a_{0})^{-1} \\
                                          &= (x+a_{0}-e^{\perp})(x+a_{0}+e^{\perp})^{-1}.
\end{align*}                                          
Hence $u \boxplus a=e+(x+a_{0}-e^{\perp})~_{e}(x+a_{0}+e^{\perp})^{-1}$. In other words, $i(e,x)\boxplus a=i(e,x+P(e^{\perp})a)$. This completes the proof. \hfill $\Box$                                          
 
 The following is now an immediate consequence of Lemma \ref{explicit} and the fact that the map $i:Y \to X$ is a bijection.                                          
  \begin{crlre}
   The map $X \ni u \to u \boxplus a \in X$ is injective for every $a \in \Omega$.
    \end{crlre}
  
  We denote the action on $Q$ on $X$, obtained by applying Lemma \ref{extension}, extending the action $\boxplus$ of $\Omega$ on $X$ by $\boxplus$ itself.
  We now prove our main theorem.
 
 \begin{thm}
 \label{main theorem}
 Let $V$ be a finite dimensional real Euclidean Jordan algebra with identity $1$ and let $Q$ be the closed convex cone of squares in $V$. Let $X:=\{u \in V: -1 \leq u \leq 1\}$. 
 Then $X$ is the Wiener-Hopf compactification of $Q$ where the action of $\Omega:=Int(Q)$ on $X$ is given by the equation
 \[
 u \boxplus a:=1-2\widetilde{u}+2P(\widetilde{u})(\widetilde{u}+a^{-1})^{-1}
 \]
 for $u \in X$ and $a \in \Omega$. Here for $u \in X$, $\widetilde{u}:=\frac{1-u}{2}$. 
 
 \end{thm}               
 \textit{Proof.} We only need to verify conditions (C1), (C2) and (C3) of Definition \ref{Wiener-Hopf compactification}.  Let $a \in Q$ be given. 
 We claim that 
 \begin{equation}
 \label{openness}
 \{u \boxplus (b +a): u \in X, b \in \Omega\}=\{u\in X: u >i(a)\}.
 \end{equation}  Note that the map $i:Q \to X$ defined by $i(x)=\frac{x-1}{x+1}$ 
 preserves the order $<$ i.e. $i(x)<i(y)$ if and only if $y-x \in \Omega$.  Let $u \in X$ and $b \in \Omega$ be given. Choose a sequence $x_{n} \in Q$
 such that $i(x_n) \to u$.   Then $i(x_{n}+b+a) \to u \boxplus (b+a)$ and $i(x_{n}+b+a) \geq i(b+a)$. Hence $u \boxplus (b+a) \geq i(b+a) > i(a)$. 
 
Let $u \in X$ be such that $u >i(a)$. Observe that there exists $\epsilon>0$ such that $u>i(a+\epsilon)$. Let $x_{n}$ be a sequence in $Q$ such that
 $i(x_{n}) \to u$. Then $i(x_{n})>i(a+\epsilon)$ eventually. Let $y_{n}:=x_{n}-a-\epsilon$. Note that $y_{n} \in \Omega$ eventually. Since $X$ is compact,
 by passing to a subsequence if necessary, we can assume that $i(y_{n})$ converges and let $v \in X$ be its limit. Then clearly \[\displaystyle v \boxplus (\epsilon +a)=\lim_{n \to \infty}i(y_{n})\boxplus (\epsilon +a)=\lim_{n\to \infty}i(y_{n}+\epsilon +a)= \lim_{n \to \infty} i(x_{n})=u.\] 
 This proves the claim. Since $\{ u \in X: u > i(a)\}$ is open in $X$, it follows that (C1) is satisfied.
 
For $u \in X$, let \[A_{u}:=\{a \in V: \textrm{ $ \exists (a_{1},a_{2},v) \in \Omega \times \Omega \times X$ such that $a=a_{1}-a_2$ and $u\boxplus a_1=v \boxplus a_2$}\}.\] 
Let $(e,x) \in Y$ be given. We claim that \[A_{i(e,x)}=\{a \in V: x+P(e^{\perp})a \geq 0\}.\]  Let $a \in A_{i(e,x)}$ be given. Then there exists $(f,y) \in Y$, $a_{1},a_{2} \in \Omega$ such that
$a=a_{1}-a_2$ and $i(e,x)\boxplus a_{1}=i(f,y)\boxplus a_{2}$. By Lemma \ref{explicit}, this implies that $e=f$ and $x+P(e^{\perp})a_{1}=y+P(e^{\perp})a_{2}$. Since $y \geq 0$, it follows that
$x+P(e^{\perp})a=x+P(e^{\perp})(a_{1})-P(e^{\perp})a_{2}=y \geq 0$.        
Let  $a \in V$ be such that $x+P(e^{\perp})a \geq 0$. Write $a=a_{1}-a_{2}$ with $a_{1},a_{2} \in \Omega$. Set $y:=x+P(e^{\perp})a$. Then $y \in V_{1}(e^{\perp}) \cap Q$.  Again, by Lemma 
\ref{explicit}, $i(e,x) \boxplus a_{1} = i(e,y) \boxplus a_{2}$. Hence $a \in A_{i(e,x)}$. This proves the claim. 

Let $u_{0}=i(0,0)$. Then $\{u_{0} \boxplus a : a \in \Omega\}=\{i(a): a \in \Omega\}$ is dense in $X$. For $i(Q)$ is dense in $X$ and $\Omega$ is dense in $Q$. 
Also note that $A_{u_{0}}=\{a \in V: 0+P(1)a \geq 0\}=Q$. This proves that (C2) is satisfied.

Let $(e,x), (f,y) \in Y$ be such that $A_{i(e,x)}=A_{i(f,y)}$. Consider an element $a$ in the subspace $V_{1}(e) \boxplus V_{\frac{1}{2}}(e)$. Note that for $x+tP(e^{\perp})a =x$ for every $t \in \mathbb{R}$. 
Hence for every $t \in \mathbb{R}$, $ta \in A_{i(e,x)}=A_{i(f,y)}$.  This implies that $y+tP(f^{\perp})a \geq 0$ for every $ t \in \mathbb{R}$. This forces that $P(f^{\perp})a=0$ i.e. $a \in V_{1}(f) \oplus V_{\frac{1}{2}}(f)$. 
That is $V_{1}(e) \oplus V_{\frac{1}{2}}(e) \subset V_{1}(f) \oplus V_{\frac{1}{2}}(f)$. By taking orthogonal complements, we obtain $V_{0}(f) \subset V_{0}(e)$.  Arguing as before with $e$ replaced by $f$, we obtain $V_{0}(e) \subset V_{0}(f)$. It follows that
$V_{1}(e^{\perp})=V_{1}(f^{\perp})$ and hence $e^{\perp}=f^{\perp}$ or $e=f$.  Now note that $-x \in A_{i(e,x)}$. Hence $-x \in A_{i(f,y)}$. This implies that $y-P(f^{\perp})x=y-P(e^{\perp})x=y-x \geq 0$. Similarly, $x-y \geq 0$. Thus it follows that $x=y$. 
Hence $(e,x)=(f,y)$. This proves that (C3) is satisfied. This completes the proof. \hfill $\Box$

\begin{rmrk}
Let $X_{0}:=\{u \boxplus a: a \in \Omega\}$. Note that  $X_{0}=\{u \in X: u >-1\}$ (Eq. \ref{openness}). Also observe that for $u \in X$, $u \in X_0$ if and only $-1$ is not in the spectrum of $u$.
\end{rmrk}

\begin{crlre}
\label{vanishing of K-theory}
Let $\mathcal{W}(Q)$ be the Wiener-Hopf $C^{*}$-algebra associated to the cone $Q$. Then the $K$-theory of $\mathcal{W}(Q)$ vanishes.
\end{crlre}
\textit{Proof.} Let $X_{0}:=\{u \boxplus a: a \in \Omega\}=\{u \in X: u >-1\}$. By Proposition 6.6 of \cite{Jean_Sundar}, it is enough to show that
the $K$-groups of $C_{0}(X_{0})$ vanish.  Note that \[C_{0}(X_{0}):=\{f \in C(X): f(u)=0 ~\textrm{if $-1\in \sigma(u)$} \}.\]For $t \in [0,1]$, let $\alpha_{t}:C(X) \to C(X)$ be 
defined by 
\[
\alpha_{t}(f)(u)=f(tu+t-1)
\]
for $f \in C(X)$ and $u \in X$. Then clearly $\alpha_{t}$ is a homotopy of $*$-homomorphisms and $\alpha_{0}(f)=f(-1)$ and $\alpha_{1}=id$. Note that $\alpha_{t}$ leaves $C_{0}(X_0)$ 
invariant. Thus $\{\alpha_{t}\}_{t \in [0,1]}$ restricted to $C_{0}(X_0)$ is a homotopy of $*$-homomorphisms connecting the identity map with the zero map. This implies that the $K$-theory
of $C_{0}(X_{0})$ vanish. This completes the proof. \hfill $\Box$

We end this paper by remarking how our approach is applicable to other examples. We consider the continuous $ax+b$-semigroup considered in \cite{Jean_Sundar} .   We refer the reader to 
\cite{Jean_Sundar} for the definition of the Wiener-Hopf $C^{*}$-algebra associated to a closed Ore subsemigroup of a locally compact group.  Let   \[G:=\Big \{\begin{pmatrix}
                                                   a & b \\
                                                   0 & 1
                                                  \end{pmatrix}: a > 0 , b \in \mathbb{R} \Big\}.\] 
 Note that $G$ is isomorphic to the semi-direct product $\mathbb{R} \rtimes (0,\infty)$ where the multiplicative group $(0,\infty)$ acts on the additive group $\mathbb{R}$ by multiplication.  Let $P$ denote the semigroup
 $[0,\infty) \rtimes [1,\infty)$ i.e. \[P:= \Big\{\begin{pmatrix} 
                                                            a & b \\
                                                            0 & 1 \end{pmatrix}: a \geq 1, b \geq 0\Big\}.\] 
Then $P$ is right Ore and solid. Let us recall the  Wiener-Hopf compactification of $P$ from \cite{Jean_Sundar}. 

 Let $Y:=[-\infty,\infty) \times [0,\infty]$. The right action of  $G$ on $Y$  is given by the formula \[
                                 (x,y).\begin{pmatrix}
        a & b \\
        0 & 1                                                                                                                
     \end{pmatrix}=\Big(\frac{x-b}{a},\frac{y}{a}\Big).
                                \]
                                
Let $X:=[-\infty,0] \times [0,1]$. Note that $X$ is $P$-invariant.  Let $X_{0}:=X(Int P)$. Observe that $X_{0}= [-\infty, 0) \times [0, 1)$. We refer the reader to Section 7 of \cite{Jean_Sundar} for the proof of the fact that the compact space
$X$, together with the above action of $P$, is the Wiener-Hopf compactification of the semigroup $P$.  

Let $\mathcal{W}(P)$ be the Wiener-Hopf $C^{*}$-algebra associated to the semigroup $P$.  We claim that the $K$-groups of $\mathcal{W}(P)$ are trivial.  It is shown in Section 7 of \cite{Jean_Sundar} that $\mathcal{W}(P)$ is Morita-equivalent
to $C_{0}(Y) \rtimes (\mathbb{R} \rtimes (0,\infty)) \cong (C_{0}(Y) \rtimes \mathbb{R})\rtimes (0,\infty)$. By Connes-Thom isomorphism, it follows that $K_{i}(\mathcal{W}(P))$ is isomorphic to $K_{i}(C_{0}(Y))$ for $i=0,1$.  
We leave it to the reader to verify that $C_{0}(Y)$ is contractible.  Thus $K_{i}(\mathcal{W}(P))=0$ for $i=0,1$. 

\begin{rmrk}
Alternatively, we could also argue as in Corollary \ref{vanishing of K-theory}. By Theorem 8.1 of \cite{Sundar_Toeplitz}, it follows that $K_{i}(\mathcal{W}(P))$ is isomorphic to $K_{i}(C_{0}(X_0))$. 
It is clear that $C_{0}(X_{0})$ is contractible. 

\end{rmrk}

\bibliography{references}

\def\cprime{$'$} \def\cprime{$'$} \def\cprime{$'$}
\providecommand{\bysame}{\leavevmode\hbox to3em{\hrulefill}\thinspace}
\providecommand{\MR}{\relax\ifhmode\unskip\space\fi MR }
\providecommand{\MRhref}[2]{%
  \href{http://www.ams.org/mathscinet-getitem?mr=#1}{#2}
}
\providecommand{\href}[2]{#2}
\begin{thebibliography}{Sun15}

\bibitem[AJ07]{AJ07}
Alexander Alldridge and Troels~Roussau Johansen, \emph{Spectrum and analytical
  indices of the {$C^\ast$}-algebra of {W}iener-{H}opf operators}, J. Funct.
  Anal. \textbf{249} (2007), no.~2, 425--453.

\bibitem[AJ08]{AJ08}
\bysame, \emph{An index theorem for {W}iener-{H}opf operators}, Adv. Math.
  \textbf{218} (2008), no.~1, 163--201.

\bibitem[All11]{All11}
Alexander Alldridge, \emph{Convex polytopes and the index of {W}iener-{H}opf
  operators}, J. Operator Theory \textbf{65} (2011), no.~1, 145--155.

\bibitem[FK94]{Faraut}
Jacques Faraut and Adam Kor{\'a}nyi, \emph{Analysis on symmetric cones}, Oxford
  Mathematical Monographs, The Clarendon Press, Oxford University Press, New
  York, 1994, Oxford Science Publications.

\bibitem[HN95]{Hilgert_Neeb}
Joachim Hilgert and Karl-Hermann Neeb, \emph{Wiener-{H}opf operators on ordered
  homogeneous spaces. {I}}, J. Funct. Anal. \textbf{132} (1995), no.~1,
  86--118.

\bibitem[MR82]{Renault_Muhly}
Paul~S. Muhly and Jean~N. Renault, \emph{{$C^{\ast} $}-algebras of
  multivariable {W}iener-{H}opf operators}, Trans. Amer. Math. Soc.
  \textbf{274} (1982), no.~1, 1--44.

\bibitem[Nic87]{Nica_WienerHopf}
Alexandru Nica, \emph{Some remarks on the groupoid approach to {W}iener-{H}opf
  operators}, J. Operator Theory \textbf{18} (1987), no.~1, 163--198.

\bibitem[RS15]{Jean_Sundar}
J~Renault and S~Sundar, \emph{Groupoids associated to {O}re semigroup actions},
  to appear in J. Operator Theory \textbf{73} (2015), no.~2, 491--514.

\bibitem[Sun15]{Sundar_Toeplitz}
S.~Sundar, \emph{Toeplitz {$C^{*}$}-algebras associated to endomorphisms of
  {O}re semigroups}, arXiv:1503.00828/math.OA, 2015.

\end{thebibliography}
\bibliographystyle{amsalpha}

\noindent
{\sc S. Sundar}
(\texttt{sundarsobers@gmail.com})\\
         {\footnotesize  Chennai Mathematical Institute, H1 Sipcot IT Park, \\
Siruseri, Padur, 603103, Tamilnadu, INDIA.}

\end{document}